\documentclass{elsart}
\usepackage{latexsym}
\usepackage{amssymb}

\def\F{\mathbb F}
\def\Z{\mathbb Z}
\def\C{\mathbb C}

\def\Tr{\mathrm{Tr}}

\begin{document}

\begin{frontmatter}
\title{The divisibility modulo 24 of Kloosterman sums on $GF(2^m)$, $m$ even}
\author{Marko Moisio}
\address{Department of Mathematics and Statistics,
Faculty of Technology, University of Vaasa, PO. Box 700, FIN-65101, Finland}
\ead{mamo@uwasa.fi}
\begin{abstract}
In a recent work by Charpin, Helleseth, and Zinoviev Kloosterman sums $K(a)$ over a finite field $\F_{2^m}$ were evaluated modulo 24 in the case $m$ odd, and the number of those $a$ giving the same value for $K(a)$ modulo 24 was given. In this paper the same is done in the case $m$ even. The key techniques used in this paper are different from those used in the aforementioned work. In particular, we exploit recent results on the number of irreducible polynomials with prescribed coefficients.    
\end{abstract} 
\end{frontmatter}

\section{Introduction}Let $\F_{2^m}$ denote the finite field of $2^m$ elements, $m>2$, and let $a\in\F_{2^m}$, $a\ne0$. Kloosterman sums $K(a)$ over $\F_{2^m}$ are widely studied for a long time for their own sake as interesting mathematical objects as well as for their connection to coding theory, most notably to the weight distribution of the Melas codes (see e.g.~\cite{La-Wo, dfhr, Moisio07} and the bibliography in them). 

The value set of $K(a)$ was obtained by Lachaud and Wolfmann~\cite{La-Wo}, and moreover, they gave the number of preimages of a specific value in terms of Kronecker class numbers. However, very little is known of the value $K(a)$ for a specific element $a$. A recent result towards to the solution of this very difficult but important problem was obtained by Charpin, Helleseth, and Zinoviev~\cite{chz} who gave congruences modulo 24 for $K(a)$ in the case $m$ odd. In this paper congruences modulo 24 are derived in the case $m$ even. The tools used in this paper are different from those used in~\cite{chz}: there $K(a)$ is linked to the number of words in a coset of weight 4 in BCH-code with minimum distance 8, and cubic exponential sums evaluated by Carlitz, but in this paper $K(a)$ is linked to the number of irreducible cubic polynomials with prescribed norm and trace, and to the number of solutions of $x^4+x^3=a$ in $\F_{2^m}$. In the calculation of the number those elements $a\in\F_q^*$ which yield the same value for $K(a)$ modulo 24 also explicit evaluations of certain exponential sums are needed, and we shall see that the value distribution of $K(a)$ modulo 24 depends on the residue class of $m$ modulo 24.

The rest of the paper is organized as follows. In Section 2 notations are fixed and some previous divisibility results needed are recalled. In Section 3 congruences modulo 3 for $K(a^4+a)$ are obtained by considering a family of elliptic curves related to the number of irreducible cubic polynomial with prescribed trace and norm. In Section 4 the number of solutions of $x^{2^k}+x^{2^k-1}=a$ in $\F_{2^m}$ is calculated. In Section 5, $K(a)$ is evaluated modulo 24 in the case $m$ even (Theorem~\ref{t:mainthm1}), and finally, in Section 6, the number of non-zero elements $a\in\F_{2^m}$ which yield the same value for $K(a)$ modulo 24 is given (Theorems~\ref{t:mainthm2} and~\ref{t:mainthm3}). 
 
\section{Preliminaries} In this section we fix some notations and recall some previous divisibility results of Kloosterman sums needed in the sequel.

Let $m>2$ be an integer and let $s$ be a positive factor of $m$. Let $q=2^m$ and let $\F_q$ denote the finite field of $q$ elements and let $\F_q^*=\F_q\setminus\{0\}$. Let $\Tr_{k}$ denote the trace function from $\F_q$ onto $\F_{2^s}$ i.e.
\[
	\Tr_s(x)=x+x^{2^s}+x^{2^{2s}}+\cdots+x^{2^{(\frac ms-1)s}}\qquad\forall\, x\in\F_q.
\]
Moreover, we use the notation $\Tr$ instead of $\Tr_1$.

Let $a\in\F_q$ and let $\chi$ be the canonical additive character of $\F_q$ i.e. $\chi(x)=(-1)^{\Tr(x)}$ for all $x\in\F_q$. We shall use frequently the following well-known result, the orthogonality of characters:
\[
	\sum_{x\in\F_q}\chi(ax)=\left\{ 
\begin{array}{ll}
	q &\quad \textrm{if $a=0$,}\\
	0 &\quad \textrm{if $a\ne0$.}
\end{array} \right.
\]

Let $K(a)$ denote the Kloosterman sum defined by
\[
K(a)=\sum_{x\in\F_q^*}\chi(x+ax^{-1}).
\]

We have the following result by Helleseth and Zinoviev:

\begin{prop}[\cite{He-Zi}]\label{p:div_mod_8}Let $a\in\F_q^*$. Then
\[
K(a)\equiv\left\{ 
\begin{array}{ll}
	3\pmod8 &\quad \textrm{if $\Tr(a)=1$,}\\
	-1\pmod8 &\quad \textrm{if $\Tr(a)=0$.}
\end{array} \right.
\]
\end{prop}

Moreover, we have the following result by Moisio:

\begin{prop}[\cite{Moisio08}]\label{p:div_mod_3}Let $a\in\F_q^*$. Then  $K(a)\equiv0\pmod3$ if and only if one of the following condition holds 
\begin{enumerate}
	\item $m$ is odd and $\Tr(\sqrt[3]{a})=0$,
	\item $m$ is even, $a=b^3$ for some $b\in\F_q$, and $\Tr_2(b)\ne0$.
\end{enumerate}
\end{prop}

\section{Kloosterman sums and irreducible cubic polynomials}
  
Let $a,b\in\F_q^*$, and let $P_3(a,b)$ denote the number of irreducible polynomials $x^3+ax^2+dx+b\in\F_q[x]$. Let $\mathcal X$ be the projective elliptic curve over $\F_q$ defined by
\[
	\mathcal X : y^2+cy+xy=x^3,
\]
where $c=b/a^3$, and let $|\mathcal X(\F_q)|$ denote the number of rational points on $\mathcal X$.  

We have the following special case of~\cite[Thm.7.3]{Moisio08}

\begin{prop}\label{p:P_3}
\[
	P_3(a,b)=\frac13(|\mathcal X(\F_q)|-\epsilon),
\]
where $\epsilon$ equals 1 or 0 depending on whether $c=1$ or $c\ne1$.
\end{prop}

By using Proposition~\ref{p:P_3} we are able to prove the main result of this section:
 
\begin{thm}\label{t:main1}Let $c\in\F_q^*$, $c\ne1$. Then
\[
K(c^4+c^3)\equiv\left\{ 
\begin{array}{ll}
	1\pmod3 &\quad \textrm{if $m$ is even and $\Tr(c)=0$,}\\
  -1\pmod3 &\quad \textrm{if $m$ is even and $\Tr(c)=1$,}\\
	0\pmod3 &\quad \textrm{if $m$ is odd.}
\end{array} \right.
\]
\end{thm}

\begin{pf}By Proposition~\ref{p:P_3} 
\[
	3P_3(1,c)=|\mathcal X(\F_q)|
\]
where $\mathcal X : y^2+cy+xy=x^3$. Write the equation of $\mathcal X$ in the form 
\[
	y^2+(x+c)y=x^3.
\] 
For a fixed $x\ne c$ substitute $y\mapsto x+c$ to get the equation in the form 
\[
	y^2+y=x^3/(x+c)^2.
\] 
Hence, by the orthogonality of characters, the number of solutions $(x,y)$ with $x\ne c$ in $\F_q^2$ of $y^2+cy+xy=x^3$ is equal to 
\begin{eqnarray*}
	N&:=&\sum_{x\in\F_q,\, x\ne c}\left(1+\chi\Bigg(\frac{x^3}{(x+c)^2}\Bigg)\right)\\
	&\stackrel{x\mapsto x+c}{=}& q-1 +\sum_{x\in\F_q^*}\chi\Bigg(\frac{(x+c)^3}{x^2}\Bigg)\\
	&=&q-1+\sum_{x\in\F_q^*}\chi(x+c+c^2x^{-1}+c^3x^{-2})\\
	&=&q-1+\chi(c)\sum_{x\in\F_q^*}\chi(x)\chi(c^2x^{-1})\chi(c^3x^{-2}).
\end{eqnarray*}
Since $\Tr(z)=\Tr(z^2)$ for all $z\in\F_q$, we now get
\[
 N=q-1+\chi(c)\sum_{x\in\F_q^*}\chi(x^2+c^4x^{-2}+c^3x^{-2})=q-1+\chi(c)K(c^4+c^3),
\]
where the last equation follows by noting that $x\mapsto x^2$ is a permutation of $\F_q$. 

Since equation $y^2+cy+xy=x^3$ has exactly one solution with $x=c$, and since $\mathcal X$ has exactly one point at infinity, we now get 
\[
	3P(1,c)=|\mathcal X(\F_q)|=q+1+\chi(c)K(c^4+c^3),
\]
and consequently $\chi(c)K(c^4+c^3)\equiv -q-1\pmod 3$. Since $q\equiv(-1)^m\pmod 3$, we get
\[
\chi(c)K(c^4+c^3)\equiv\left\{ 
\begin{array}{ll}
	1\pmod3 &\quad \textrm{if $m$ is even,}\\
	0\pmod3 &\quad \textrm{if $m$ is odd.}
\end{array} \right.
\]
This completes the proof, since $\chi(c)=1,-1$ depending on whether $\Tr(c)=0,1$. \qed
\end{pf}

\begin{rem} If we can prove in the case $m$ even that all such elements $a$ in $\F_q^*$, which do not satisfy the condition in Proposition~\ref{p:div_mod_3}, can be written in the form $a=c^4+c^3$, then Proposition~\ref{p:div_mod_3} and Theorem~\ref{t:main1} give the divisibility modulo 3 of $K(a)$ for all $\in\F_q^*$. We shall see in the next section that this indeed is the case.
\end{rem}

\section{The equation $x^{2^k}+x^{2^k-1}=a$}

In this section we shall prove the following 

\begin{thm}\label{t:main2}Let $k$ be a positive integer and let $a\in\F_q^*$. The number $N(a)$ of solutions of
\begin{equation}\label{e:key}
	x^{2^k}+x^{2^k-1}=a
\end{equation}
is given by
\[
	N(a)=\left\{
\begin{array}{ll}
	1 &\quad \textrm{if $a\ne b^{2^k-1}$ for all $b\in\F_q$,}\\
	2^s &\quad \textrm{if $a=b^{2^k-1}$ for some $b\in\F_q$, and $\Tr_s(b)=0$,}\\
	0 &\quad \textrm{if $a=b^{2^k-1}$ for some $b\in\F_q$, and $\Tr_s(b)\ne0$,}
\end{array} \right.	
\]
where $s=\gcd(k,m)$.
\end{thm}

\begin{pf}
Substitute $x\mapsto x^{-1}$ to~(\ref{e:key}) and then multiply both sides by $x^{2^k}$ to get an equivalent equation
\[
	1+x=ax^{2^k}. 
\]
Now, by the orthogonality of characters, we get
\begin{eqnarray*}
	qN(a)&=&\sum_{c\in\F_q}\sum_{x\in\F_q}\chi(c(1+x+ax^{2^k}))\\
	&=&q+\sum_{c\in\F_q^*}\chi(c)\sum_{x\in\F_q}\chi(cx)\chi(cax^{2^k}). 
\end{eqnarray*}
   
Since $\chi(z^{2^k})=\chi(z)$ for all $z\in\F_q$, we obtain
\[
	qN(a)=q+\sum_{c\in\F_q^*}\chi(c)\sum_{x\in\F_q}\chi((c^{2^k}+ca)x^{2^k}). 
\]

Since $x\mapsto x^{2^k}$ is a permutation of $\F_q$, the orthogonality of characters implies that the inner sum
\[
	\sum_{x\in\F_q}\chi((c^{2^k}+ca)x^{2^k})=\left\{
\begin{array}{ll}
	0 &\quad \textrm{if $c^{2^k-1}\ne a$,}\\
	q &\quad \textrm{if $c^{2^k-1}=a$.} 
\end{array} \right.	
\]
  
Hence, if $c^{2^k-1}=a$ is not solvable, then $N(a)=1$. On the other hand, $c^{2^k-1}=a$ is solvable if and only if $a$ is in $<\gamma^{2^s-1}>$ since $\gcd(2^k-1,2^m-1)=2^s-1$. Moreover, if $b$ is one solution, then all the solutions are $b\alpha$ where $\alpha$ runs over $\F_{2^s}^*$. Hence, if $a=b^{2^k-1}$ for some $b\in\F_q^*$, then
\[
	qN(a)=q+q\sum_{\alpha\in\F_{2^s}^*}\chi(b\alpha)=\left\{
\begin{array}{ll}
	q2^s&\quad \textrm{if $\Tr_s(b)=0$,}\\
	0 &\quad \textrm{if $\Tr_s(b)\ne0$,} 
\end{array} \right.
\]
and therefore, in this case,
\[
	N(a)=\left\{
\begin{array}{ll}
	2^s&\quad \textrm{if $\Tr_s(b)=0$,}\\
	0&\quad \textrm{if $\Tr_s(b)\ne0$,} 
\end{array} \right.
\]
completing the proof.\qed
\end{pf}

\begin{cor}\label{c:N_odd}Assume $\gcd(k,m)=1$. Then 
\[
	N(a)=\left\{
\begin{array}{ll}
	2 &\quad \textrm{if $\Tr(a^{\frac1{2^k-1}})=0$,}\\
	0 &\quad \textrm{otherwise.}
\end{array} \right.	
\]
\end{cor}

\begin{rem}Corollary~\ref{c:N_odd} was also proved in~\cite{Gar-Lis} by using different methods. Moreover, Corollary~\ref{c:N_odd} and Proposition~\ref{p:div_mod_3} imply that $K(a)$ is divisible by 3 if and only if 
$a = c^4 + c^3$ for some $c\in\F_q^*$. This result was also proved in~\cite{Gar-Lis}.
\end{rem}

\section{The evaluation of $K(a)$ modulo 24, $m$ even}

We are now able to evaluate $K(a)$ modulo 24:

\begin{thm}\label{t:mainthm1} Let $a\in\F_q^*$.
\begin{enumerate}
\item Assume $a=b^3$ for some $b\in\F_q$ with $\Tr_2(b)\ne0$. Then, $a \ne c^4+c^3$ for all $c\in\F_q$, and
\[
K(a)\equiv\left\{ 
\begin{array}{ll}
	15\pmod{24} &\quad \textrm{if $\Tr(a)=0$,}\\
	3\pmod{24} &\quad \textrm{if $\Tr(a)=1$.}
\end{array} \right.
\]
\item Otherwise, $a=c^4+c^3$ for some $c\in\F_q$, and
\[
K(a)\equiv\left\{ 
\begin{array}{ll}
  7\pmod{24} &\quad \textrm{if $\Tr(c)=0$ and $\Tr(c^3)=0$,}\\
  19\pmod{24} &\quad \textrm{if $\Tr(c)=0$ and $\Tr(c^3)=1$,}\\
  11\pmod{24} &\quad \textrm{if $\Tr(c)=1$ and $\Tr(c^3)=0$,}\\ 
  23\pmod{24} &\quad \textrm{if $\Tr(c)=1$ and $\Tr(c^3)=1$.} 
\end{array} \right.
\]	
\end{enumerate}
\end{thm} 

\begin{pf} If $a=b^3$ for some $b\in\F_q$ and $\Tr_2(b)=0$ or $a\ne b^3$ for all $b\in\F_q$, then $a=c^4+c^3$ for some $c\in\F_q$ by Theorem~\ref{t:main2} ($k=2$). Now, by Theorem~\ref{t:main1}, 
\[
K(a)\equiv\left\{ 
\begin{array}{ll}
	1\pmod{3} &\quad \textrm{if $\Tr(c)=0$,}\\
	-1\pmod{3} &\quad \textrm{if $\Tr(c)=1$.}
\end{array} \right.
\]
These congruences and the congruences in Proposition~\ref{p:div_mod_8} imply
\[
K(a)\equiv\left\{ 
\begin{array}{ll}
	7\pmod{24} &\quad \textrm{if $\Tr(a)=0$ and $\Tr(c)=0$,}\\
	23\pmod{24} &\quad \textrm{if $\Tr(a)=0$ and $\Tr(c)=1$,}\\
  19\pmod{24} &\quad \textrm{if $\Tr(a)=1$ and $\Tr(c)=0$,}\\ 
  11\pmod{24} &\quad \textrm{if $\Tr(a)=1$ and $\Tr(c)=1$.}
\end{array} \right.
\]
Since $\Tr(a)=\Tr(c^4)+\Tr(c^3)=\Tr(c)+\Tr(c^3)$, the proof is complete in this case. 

In the remaining case we have $a \ne c^4+c^3$ for all $c\in\F_q$, by Theorem~\ref{t:main2}, and the congruences in Propositions~\ref{p:div_mod_8} and~\ref{p:div_mod_3} complete the proof. \qed    
\end{pf}

\section{The value distribution of $K(a)$ modulo 24, $m$ even}

Assume that $m$ is even, and let $a\in\F_q^*$. In this section we shall give the cardinality of the preimage of the remainder $K(a)$ mod 24 under the map $y\mapsto (K(y)\ \mathrm{ mod }\ 24)$ defined on $\F_q^*$. We split our consideration according to whether $a$ can, or can not, be represented in the form $a=c^4+c^3$ for some $c\in\F_q.$

We shall need the following explicit result on exponential sums.

\begin{lem}\label{l:expsums}
\begin{eqnarray*} 
&\sum_{x\in\F_q}\chi(x^3)&=-(-1)^{\frac m2}2\sqrt q,\\
&\sum_{x\in\F_q}\chi(x^9)&=\left\{ 
\begin{array}{ll}
	-(-1)^{\frac m2}2\sqrt q&\quad \textrm{if $m\equiv\pm1\pmod 3$,}\\
	-(-1)^{\frac m2}8\sqrt q&\quad \textrm{if $m\equiv 0\pmod 3$,}
\end{array} \right.\\
&\sum_{x\in\F_q}\chi(x^3+x)&=\left\{
\begin{array}{ll}
	-2\sqrt q&\quad \textrm{if $m\equiv0\pmod 8$,}\\
	0&\quad \textrm{if $m\equiv2,6\pmod 8$,}\\
	2\sqrt q&\quad \textrm{if $m\equiv4\pmod 8$.}
\end{array} \right.\\
&\sum_{x\in\F_q}\chi(x^9+x^3)&=\left\{
\begin{array}{ll}
	-8\sqrt q&\quad \textrm{if $m\equiv0\pmod 8$,}\\
	2\sqrt q&\quad \textrm{if $m\equiv2,6\pmod 8$,}\\
	4\sqrt q&\quad \textrm{if $m\equiv4\pmod 8$.}
\end{array} \right.
\end{eqnarray*}
\end{lem}

\begin{pf}Let $S(f)=\sum_{x\in\F_q^*}\chi(f(x))$. If $f(x)=x^3$, or $f(x)=x^9$ and $m\equiv 0\pmod 3$, the result is well-known (see e.g~\cite{Wolfmann, Moisio00}). If $m\equiv\pm1\pmod3$ then $\gcd(9,q-1)=3$, and therefore $S(x^9)=S(x^3)$. If $f(x)=x^3+x$, we refer to~\cite[p.191]{Stichtenoth}. 

Assume $f(x)=x^9+x^3$ and consider the projective Artin-Schreier curve $\mathcal X$ over $\F_2$ defined by 
$\mathcal X:y^2+y=x^9+x^3$. The genus of $\mathcal X$ is $g=4$, and the number of $\F_{2^r}$-rational points 
$|\mathcal X(\F_{2^r})|$ on $\mathcal X$ is given by
\[
	|\mathcal X(\F_{2^r})|=2^r+1+\underbrace{\sum_{x\in\F_{2^r}}(x^9+x^3)}_{:=S_r(f)}.
\]

On the other hand,
\[
	|\mathcal X(\F_{2^r})|=2^r+1-\sum_{i=1}^8\omega_i^{-r},
\] 
where the $\omega_i$ are the roots of the $L$-polynomial $L(t)=a_0+a_1+\cdots+a_8t^8$ of $\mathcal X$ in $\C$ (see ~\cite[p.168]{Stichtenoth}). Moreover, the coefficients $a_i$ satisfy the recursion: $a_0=1$ and
\[
	ia_i=S_i(f)a_0+S_{i-1}(f)a_1+\cdots+S_1(f)\qquad\mathrm{for }\ i=1,2,3,4.
\]   
We also know that $a_{8-i}=2^{4-i}a_i$ for $i=0,1,2,3$. 

It is easy to see that $S_1(f)=2, S_2(f)=4, S_3(f)=-4$, and $S_4(f)=16$. It now follows that
\begin{eqnarray*}
L(t)&=&16t^8+16t^7+16t^6+8t^5+8t^4+4t^3+4t^2+2t+1\\
&=&(2t^2-2t+1)(2t^2+2t+1)^2(2t^2+1).
\end{eqnarray*}
The reciprocals of the roots of the factors of $L(t)$ are $1\pm\sqrt{-1}$, $-1\pm\sqrt{-1}$, and $\pm\sqrt{-2}$, respectively, and it follows that
\begin{eqnarray*}
	\sum_{x\in\F_{2^m}}(x^9+x^3)&=-&\textstyle\Big(2\cdot2^{\frac m2}\cos(\frac{m\pi}4)
	+4\cdot(-1)^m2^{\frac m2}\cos(\frac{m\pi}4)+2\cdot2^{\frac m2}\cos(\frac{m\pi}2)\Big)\\
	&=&-\textstyle2\cdot2^{\frac m2}\Big(3\cos(\frac{m\pi}4)+\cos(\frac{m\pi}2)\Big).	
\end{eqnarray*}
The claimed formula follows now easily.\qed
\end{pf}

\subsection{Case $a=c^4+c^3$ for some $c\in\F_q$}

Let $\gamma$ be a primitive element of $\F_q$. Let $\epsilon,\delta\in\F_2$, and let 
\[
	C(\epsilon,\delta)=\{c\in\F_q^*\setminus\{1\}\mid \Tr(c)=\epsilon,\, \Tr(c^3)=\delta\}.
\]
Consider the function $f$ on $C(\epsilon,\delta)$ defined by the polynomial $f(x)=x^4+x^3$. Let $N(\epsilon,\delta)$ denote the number of elements $c$ in $C(\epsilon,\delta)$ satisfying $f(c)\in<\gamma^3>$. 

\begin{lem}\label{l:N_N'}The value of $K(f(c))$ modulo 24 corresponding to the pair $(\Tr(c),\Tr(c^3))=(\epsilon,\delta)$ in Theorem~\ref{t:mainthm1} is attained exactly
\[
 \#C(\epsilon,\delta)-\textstyle\frac34N(\epsilon,\delta)
\]
times as $c$ varies over $C(\epsilon,\delta)$. Moreover, 
\[
	N(\epsilon,\delta)=\#\{i=1,\dots,\textstyle\frac{q-1}3-1 \mid \Tr(\gamma^{3i})=\epsilon,\, \Tr(\gamma^{9i})=\delta\}.
\] 
\end{lem}

\begin{pf} Let $a$ be an element in the image of $f$. Assume $a\in<\gamma^3>$. Now, by Theorem~\ref{t:main2}, there are exactly four elements $c\in\F_q^*$ such that $f(c)=a$. Each such $c$ must belong to $C(\epsilon,\delta)$, for otherwise $\Tr(c)\ne\epsilon$ or $\Tr(c^3)\ne\delta$ leading to a different value $K(a)$ mod 24, by Theorem~\ref{t:mainthm1}. 

If $a\not\in<\gamma^3>$, then by Theorem~\ref{t:main2}, $a$ has exactly one preimage under $f$ and therefore $K(a)$ mod 24 is attained exactly $\frac14N(\epsilon,\delta)+N'(\epsilon,\delta)$ times, where $N'(\epsilon,\delta)$ is the number of elements $c$ in $C(\epsilon,\delta)$ satisfying $f(c)\not\in<\gamma^3>$. But $N'(\epsilon,\delta)=\#C(\epsilon,\delta)-N(\epsilon,\delta)$, which proves the first part of the lemma. 

To prove the claimed expression for $N(\epsilon,\delta)$, we note that $a=c^3(c+1)\in<\gamma^3>$ if and only if $c+1\in<\gamma^3>$ if and only if $c=\gamma^{3i}+1$ for some $i=1,\dots,(q-1)/3-1$. If $c=\gamma^{3i}+1$, then $c^3=\gamma^{9i}+\gamma^{6i}+\gamma^{3i}+1$, and consequently 
\[
	\Tr(c)=\Tr(\gamma^{3i})+\Tr(1)=\Tr(\gamma^{3i})
\] 
and
\[
	\Tr(c^3)=\Tr(\gamma^{9i})+\Tr((\gamma^{3i})^2)+\Tr(\gamma^{3i})+\Tr(1)=\Tr(\gamma^{9i}).
\]

Hence, the number $N(\epsilon,\delta)$ of elements $c$ in $C(\epsilon,\delta)$ satisfying $f(c)\in<\gamma^3>$ is given by
\[
	N(\epsilon,\delta)=\{i=1,\dots,\textstyle\frac{q-1}3-1 \mid \Tr(\gamma^{3i})=\epsilon,\, \Tr(\gamma^{9i})=\delta\},
\]
which completes the proof.\qed
\end{pf}

Next we shall find exponential sum expressions for the numbers $N(\epsilon,\delta)$ and $\#C(\epsilon,\delta)$.
\begin{lem}\label{l:exprep}We have
\begin{eqnarray*}
	12N(\epsilon,\delta)&=&q+(-1)^{\delta}\sum_{x\in\F_q}\chi(x^9)+(-1)^{\epsilon}\sum_{x\in\F_q}\chi(x^3)
	+(-1)^{\epsilon+\delta}\sum_{x\in\F_q}\chi(x^9+x^3)-4h,
\end{eqnarray*}
and
\[	
	4\cdot\#C(\epsilon,\delta)=q+(-1)^{\delta}\sum_{x\in\F_q}\chi(x^3)+(-1)^{\epsilon+\delta}\sum_{x\in\F_q}\chi(x^3+x)-2h,
\]
where $h=4$ if $\epsilon=\delta=0$, and otherwise $h=0$.
\end{lem}

\begin{pf} Let us first calculate $N(\epsilon,\delta)$. Let $z\in\F_q$ satisfying $\Tr(z)=1$ and let $\psi$ be the canonical additive character of $\F_2$.  By the orthogonality of characters
\begin{eqnarray*}
	4N(\epsilon,\delta)&=&\sum_{i=1}^{\frac{q-1}3-1}\Big(\sum_{u\in\F_2}\psi(\Tr(\gamma^{3i}+z\epsilon)u)\Big)
	\Big(\sum_{v\in\F_2}\psi(\Tr(\gamma^{9i}+z\delta)v)\Big)\\
	&=&\sum_{i=1}^{\frac{q-1}3-1}\Big(1+(-1)^{\epsilon}\chi(\gamma^{3i})\Big)
	\Big(1+(-1)^{\delta}\chi(\gamma^{9i})\Big)\\
	&=&\sum_{i=1}^{\frac{q-1}3-1}\Big(1+(-1)^{\delta}\chi(\gamma^{9i})+(-1)^{\epsilon}\chi(\gamma^{3i})+
	(-1)^{\epsilon+\delta}\chi(\gamma^{9i}+\gamma^{3i})\Big)\\
	&=&\sum_{i=0}^{\frac{q-1}3-1}\Big(1+(-1)^{\delta}\chi(\gamma^{9i})+(-1)^{\epsilon}\chi(\gamma^{3i})+
	(-1)^{\epsilon+\delta}\chi(\gamma^{9i}+\gamma^{3i})\Big)-h,
\end{eqnarray*}
where $h=1+(-1)^{\delta}+(-1)^{\epsilon}+(-1)^{\epsilon+\delta}$. Since the values of $\gamma^{3i}$ and $\gamma^{9i}$  depend only on the residue class modulo $(q-1)/3$ of $i$, we now get
\begin{eqnarray*}
4N(\epsilon,\delta)+h&=&\frac13\sum_{x\in\F_q^*}\Big(1+(-1)^{\delta}\chi(x^9)+(-1)^{\epsilon}\chi(x^3)+
	(-1)^{\epsilon+\delta}\chi(x^9+x^3)\Big)\\
	&=&\frac13\sum_{x\in\F_q}\Big(1+(-1)^{\delta}\chi(x^9)+(-1)^{\epsilon}\chi(x^3)+
	(-1)^{\epsilon+\delta}\chi(x^9+x^3)\Big)-\textstyle\frac h3,
\end{eqnarray*}
from which the claimed formula for $N(\epsilon,\delta)$ follows.

By the orthogonality of characters we also get
\begin{eqnarray*}
	4\cdot\#C(\epsilon,\delta)&=&\sum_{i=1}^{q-2}\Big(\sum_{u\in\F_2}\psi(\Tr(\gamma^i+\alpha\epsilon)u)\Big)
	\Big(\sum_{v\in\F_2}\psi(\Tr(\gamma^{3i}+\alpha\delta)v)\Big)\\
%	&=&\sum_{x\in\F_q^*}\Big(1+(-1)^{\epsilon}\chi(x)\Big)\Big(1+(-1)^{\delta}\chi(x^3)\Big)\\
	&=&\sum_{i=1}^{q-2}\Big(1+(-1)^{\delta}\chi(\gamma^{3i})+(-1)^{\epsilon}\chi(\gamma^i)
	+(-1)^{\epsilon+\delta}\chi(\gamma^{3i}+\gamma^i)\Big)\\
	&=&\sum_{x\in\F_q}\Big(1+(-1)^{\delta}\chi(x^3)+(-1)^{\epsilon}\chi(x)+(-1)^{\epsilon+\delta}\chi(x^3+x)\Big)-2h\\
	&=&q+(-1)^{\delta}\sum_{x\in\F_q}\chi(x^3)+(-1)^{\epsilon+\delta}\sum_{x\in\F_q}\chi(x^3+x)-2h,
\end{eqnarray*}
since $\sum_{x\in\F_q}\chi(x)=0$. The proof is now complete.\qed
\end{pf}

\begin{thm}\label{t:mainthm2}Let $k\in\{7,19,11,23\}$. The number $N(k)$ of elements $a$ in $\F_q^*$ for which $K(a)\equiv k\pmod {24}$ is given by

\begin{tabular}{c|c|c|c|c}
$m\ \mathrm{mod\ 24}$ & $N(7)-3\cdot2^{m-4}$ & $N(19)-3\cdot2^{m-4}$ & $N(11)-3\cdot2^{m-4}$ & $N(23)-3\cdot2^{m-4}$\\
\hline
$0$ & $2^{\frac m2-3}-1$ & $2^{\frac m2-3}$ & $-2^{\frac m2-3}$ & $-2^{\frac m2-3}$\\
$\pm6$ & $-2^{\frac m2-2}-1$ & $0$ & $2^{\frac m2-2}$ & $0$\\
$12$ & $3\cdot2^{\frac m2-3}-1$ & $-2^{\frac m2-3}$ & $-3\cdot2^{\frac m2-3}$ & $2^{\frac m2-3}$\\
$\pm8$ & $-2^{\frac m2-2}-1$ & $2^{\frac m2-1}$ & $-2^{\frac m2-1}$ & $2^{\frac m2-2}$\\
$\pm2,\pm10$ & $2^{\frac m2-3}-1$ & $-3\cdot2^{\frac m2-3}$ & $5\cdot2^{\frac m2-3}$ & $-3\cdot2^{\frac m2-3}$\\
$\pm 4$ & $-1$ & $2^{\frac m2-2}$ & $-3\cdot2^{\frac m2-2}$ & $2^{\frac m2-1}$
\end{tabular}
\end{thm}

\begin{pf}Combine Lemmas~\ref{l:exprep} and~\ref{l:expsums} to get the following tables for $\frac34N(\epsilon,\delta)$ and $\#C(\epsilon,\delta)$:

\begin{tabular}{c|c|c|c|c}
$m$ mod 24 & $\frac34N(0,0)-2^{m-4}$ & $\frac34N(0,1)-2^{m-4}$ & $\frac34N(1,0)-2^{m-4}$ & $\frac34N(1,1)-2^{m-4}$\\
\hline
$0$ & $-9\cdot2^{\frac m2-3}-1$ & $7\cdot2^{\frac m2-3}$ & $2^{\frac m2-3}$ & $2^{\frac m2-3}$\\
$\pm6$ & $3\cdot2^{\frac m2-2}-1$ & $-2^{\frac m2-1}$ & $2^{\frac m2-2}$ & $-2^{\frac m2-1}$\\
$12$ & $-3\cdot2^{\frac m2-3}-1$ & $2^{\frac m2-3}$ & $-5\cdot2^{\frac m2-3}$ & $7\cdot2^{\frac m2-3}$\\
$\pm8$ & $-3\cdot2^{\frac m2-2}-1$ & $2^{\frac m2-1}$ & $2^{\frac m2-1}$ & $-2^{\frac m2-2}$\\
$\pm2,\pm10$ & $3\cdot2^{\frac m2-3}-1$ & $-2^{\frac m2-3}$ & $-2^{\frac m2-3}$ & $-2^{\frac m2-3}$\\
$\pm 4$ & $-1$ & $-2^{\frac m2-2}$ & $-2^{\frac m2-2}$ & $2^{\frac m2-1}$
\end{tabular} 

\begin{tabular}{c|c|c|c|c}
$m$ mod 8 & $\#C(0,0)-2^{m-2}$ & $\#C(0,1)-2^{m-2}$ & $\#C(1,0)-2^{m-2}$ & $\#C(1,1)-2^{m-2}$\\
\hline
$0$ & $-2^{\frac m2}-2$ & $2^{\frac m2}$ & $0$ & $0$\\
$\pm2$ & $2^{\frac m2-1}-2$ & $-2^{\frac m2-1}$ & $2^{\frac m2-1}$ & $-2^{\frac m2-1}$\\
$4$ & $-2$ & 0 & $-2^{\frac m2}$ & $2^{\frac m2}$
\end{tabular}

The definitions of $\epsilon$ and $\delta$ together with Lemma~\ref{l:N_N'} and Theorem~\ref{t:mainthm1} now completes the proof.\qed
\end{pf}

\begin{exmp}\label{ex:verify}Let $m=6$. By~\cite{La-Wo} we know that that the value set of $K(a)$ is $S:=\{-13,-9,-5,-1,3,7,11,15\}$. Moreover, each value $t$ in $S$ is attained exactly $H(t^2-256)$ times, where $H(d)$ is the Kronecker class number of $d$. Hence, we have the following table

\begin{tabular}{c|cccccccc}
$t$ & $-13$ & $-9$ & $-5$ & $-1$ & $3$ & $7$ & $11$ & $15$\\
\hline
$t\ \mathrm{mod}\ 24$ & $11$ & $15$ & $19$ & $23$ & $3$ & $7$ & $11$ & $15$\\
\hline
$H(t^2-256)$ & $6$ & $7$ & $12$ & $12$ & $6$ & $9$ & $8$ & $3$\\
\end{tabular} 

By Theorem~\ref{t:mainthm2}, $N(7)=3\cdot2^{6-4}-2^{\frac 62-2}-1=9$, $N(19)=3\cdot2^{2}=12$, $N(11)=3\cdot2^{2}+2^{\frac 62-2}=14=6+8$, and $N(23)=N(19)=12$. This is in accordance with the table above. The remaining values $N(3)$ and $N(15)$ will be verified in the next subsection.
\end{exmp}

\subsection{Case $a \ne c^4+c^3$ for all $c\in\F_q$}

Assume $a \ne c^4+c^3$ for all $c\in\F_q$, equivalently $a=b^3$ for some $b\in\F_q$ with $\Tr_2(b)\ne0$. Let $\epsilon\in\F_2$, let $\beta$ be an element of $\F_4^*$, and let
\[
	S_{\beta}(\epsilon)=\{b\in\F_q^*\mid \Tr(b^3)=\epsilon, \Tr_2(b)=\beta\}.
\]

\begin{lem}\label{l:N_b_cubic}The value of $K(b^3)$ modulo 24 corresponding to $\Tr(b^3)=\epsilon$ in Theorem~\ref{t:mainthm1} is attained exactly
\[
	\frac13\sum_{\beta\in\F_4^*}\#S_{\beta}(\epsilon)
\]
times as $b$ varies over $\bigcup_{\beta\in\F_4^*}S_{\beta}(\epsilon)$.
\end{lem}

\begin{pf}Let $\zeta$ be a primitive element of $\F_4^*$, and let $\beta\in\F_4^*$. If $b\in S_{\beta}(\epsilon)$ and $i\in\Z$, then $\Tr((\zeta^i b)^3)=\Tr(b^3)=\epsilon$ and $\Tr_2(\zeta^i b)=\zeta^i\Tr_2(b)$, and therefore 
$\zeta^i b\in S_{\zeta^i\beta}(\epsilon)$. Hence, if $a=b^3$ for some $b\in S_{\beta}(\epsilon)$, then $a$ has exactly three preimages under the map $x\mapsto x^3$ defined on $\bigcup_{\beta\in\F_4^*}S_{\beta}(\epsilon)$. This completes the proof. \qed   
\end{pf}

\begin{thm}\label{t:mainthm3}Let $k\in\{3,15\}$. The number $N(k)$ of elements $a$ in $\F_q^*$ for which $K(a)\equiv k\pmod {24}$ is given by

\begin{tabular}{c|c|c}
$m \pmod8$ & $N(3)$ & $N(15)$\\
\hline
$0$ & $2^{m-3}$ & $2^{m-3}$\\
$2,6$ & $2^{m-3}-2^{\frac m2-2}$ & $2^{m-3}+2^{\frac m2-2}$\\
$4$ & $2^{m-3}+2^{\frac m2-1}$ & $2^{m-3}-2^{\frac m2-1}$
\end{tabular}
\end{thm}

\begin{pf}Let $\beta\in\F_4^*$, and let $\psi$ and $\eta$ be the canonical additive characters of $\F_2$ and $\F_4$. Let $z,w\in\F_q$ satisfying $\Tr(z)=1$ and $\Tr_2(w)=1$. The orthogonality of characters implies
\begin{eqnarray*}
8\cdot\#S_{\beta}(\epsilon)&=&\sum_{x\in\F_q^*}\Big(\sum_{u\in\F_2}\psi(\Tr(x^3+z\epsilon)u)\Big)
\Big(\sum_{v\in\F_4}\psi(\Tr_2(x+w\beta)v)\Big)\\
&=&\sum_{x\in\F_q^*}\Big(1+(-1)^{\epsilon}\chi(x^3)\Big)\sum_{v\in\F_4}\chi(xv)\psi(\beta v)).	
\end{eqnarray*}
Now, since $v^3=1$ for $v\in\F_4^*$, we get
\begin{eqnarray*}
8\sum_{\beta\in\F_4^*}\#S_{\beta}(\epsilon)&=&\sum_{x\in\F_q^*}\Big(1+(-1)^{\epsilon}\chi(x^3)\Big)
\sum_{v\in\F_4}\chi(xv)\sum_{\beta\in\F_4^*}\psi(\beta v)\\
&=&\sum_{x\in\F_q^*}\Big(1+(-1)^{\epsilon}\chi(x^3)\Big)\Big(3-\sum_{v\in\F_4^*}\chi(vx)\Big)\\
&\stackrel{x\mapsto v^{-1}x}{=}&\sum_{x\in\F_q^*}\Big(1+(-1)^{\epsilon}\chi(x^3)\Big)
\Big(3-\sum_{v\in\F_4^*}\chi(x)\Big)\\
&=&3\sum_{x\in\F_q^*}\Big(1+(-1)^{\epsilon}\chi(x^3)\Big)\Big(1-\chi(x)\Big).
\end{eqnarray*} 
Since $\sum_{x\in\F_q^*}\chi(x)=-1$, it follows that
\[	8\sum_{\beta\in\F_4^*}\#S_{\beta}(\epsilon)=3\Big(q+(-1)^{\epsilon}\sum_{x\in\F_q^*}\chi(x^3)-(-1)^{\epsilon}\sum_{x\in\F_q^*}\chi(x^3+x)\Big). 
\]
Now, by Lemma~\ref{l:expsums}, we get
\[
8\sum_{\beta\in\F_4^*}\#S_{\beta}(\epsilon)=3\left\{
\begin{array}{ll}
	q&\quad \textrm{if $m\equiv0\pmod 8$,}\\
	q+(-1)^{\epsilon}2\sqrt q&\quad \textrm{if $m\equiv2,6\pmod 8$,}\\
	q-(-1)^{\epsilon}4\sqrt q&\quad \textrm{if $m\equiv4\pmod 8$.}
\end{array}\right.
\]
Lemma~\ref{l:N_b_cubic} now completes the proof.\qed
\end{pf}

\begin{exmp}Let $m=6$. By Theorem~\ref{t:mainthm3}, $N(3)=2^{6-3}-2^{\frac 62-2}=6$ and 
$N(15)=2^{6-3}+2^{\frac 62-2}=10=7+3$. This is in accordance with the table in Example~\ref{ex:verify}.
\end{exmp}

\end{document}